\documentclass[12pt,reqno,notitlepage]{amsart}
\usepackage{amsmath,amsfonts,amsthm,amssymb}
\usepackage{enumerate}
\usepackage{indentfirst}
\usepackage[dvips]{graphicx}
\usepackage[all]{xy}
\usepackage{stackrel}
\usepackage{marginnote}

\usepackage[centering, includeheadfoot, hmargin=1.0in, tmargin=1.0in, 
bmargin=1in, headheight=6pt]{geometry}
\usepackage[latin1]{inputenc}
\usepackage[T1]{fontenc}
\usepackage{verbatim}
\usepackage{color}
\usepackage[normalem]{ulem}

\usepackage{hyperref}  
\hypersetup{
pdfborder={0 0 0}, 
colorlinks=true, 
citecolor=blue,
linktoc=page,
pdfauthor={Thiago Fassarella and Frank Loray}, 
pdftitle={Higgs bundles}
}
\renewcommand{\ref}{\hyperref}

\newtheorem{thm}{Theorem}[section]

\newtheorem{prop}[thm]{Proposition}

\newtheorem{lemma}[thm]{Lemma}

\newtheorem{cor}[thm]{Corollary}
\newtheorem{remark}[thm]{Remark}

\theoremstyle{definition}

\DeclareMathOperator{\tr}{tr}

\renewcommand{\P}{\mathbb{P}}
\newcommand{\C}{\mathbb{C}}

\numberwithin{equation}{section}

\newcommand{\Bun}{\mathrm{Bun}}
\newcommand{\End}{\mathcal{E}nd}
\newcommand{\SEnd}{\mathcal{SE}nd}
\newcommand{\Hom}{\mathcal{H}om}
\newcommand{\elem}{\mathrm{elem}}

\begin{document}
\title{Hitchin fibration under ramified coverings}
\author[T. Fassarella]{Thiago Fassarella}
\address{\color{black} Universidade Federal Fluminense, Instituto de Matem\'atica e Estat\'istica\\
Rua Alexandre Moura 8, S\~ao Domingos, 24210-200 Niter\'oi RJ,
Brazil.}
\email{\color{black}tfassarella@id.uff.br}

\author[F. Loray]{Frank Loray}
\address{Univ Rennes, CNRS, IRMAR - UMR 6625, F-35000 Rennes, France}
\email{frank.loray@univ-rennes1.fr}

\thanks{The second author is supported by CNRS and Henri Lebesgue Center,
 program ANR-11-LABX-0020-0. The authors also thank  Brazilian-French Network in Mathematics and CAPES-COFECUB  project 932/19.}
 \date{\today}

\begin{abstract}
We are interested in studying the variation of the Hitchin fibration in moduli spaces of parabolic Higgs bundles, under the action of a ramified covering. Given a degree two map $\pi: Y\to X$ between compact Riemann surfaces, we may pull back a Higgs bundle from $X$ to  $Y$, the lifted Higgs bundle tends to have many apparent singularities, then we perform a suitable birational transformation in order to eliminate them. This correspondence preserves the Hitchin fibrations and then its restriction to a general fiber gives a map between Abelian varieties. The aim of this paper is to describe this map. 
\end{abstract}

\maketitle


\section{Introduction}

Let $X$ be a compact Riemann surface of genus $g_X$ and let $D = t_1+\cdots + t_n$ be a divisor on it determined by $n$ distinct points. These points are called parabolic points.  We are interested in moduli spaces of Higgs bundles on $X$ with simple poles on $D$. Higgs bundles were originally  introduced by Hitchin \cite{Hi1, Hi2}  (see also  \cite{Ni}), without parabolic structure, and then extended to parabolic context in \cite{Si, Yo93, Yo95, Ko, Na}. The Hitchin fibration is defined on these moduli spaces and it is a source of remarkable applications \cite{HT,KW,Ngo1,Ngo2}. 

Before going into technical definitions let us offer some motivation for the techniques developed in this paper. 
On the one hand, it is well known that the general fiber of the Hitchin fibration is isomorphic to an Abelian variety. On the other hand, singular fibers are difficult to deal with and  this constitutes an interesting problem in the area.  We believe that the symmetries given by elementary transformations play an important role in this study. Therefore,  the first objective of this paper is to establish a dictionary that allows us to expressing elementary transformations in terms of the Beauville-Narasimhan-Ramanan (BNR) correspondence. 
As a consequence, we present a study of the variation of the Hitchin fibration with respect to a degree two ramified covering. 
Our main motivation is the application of these techniques to explicit examples, which are very rare in the current literature. Hence, the cases $(g_X,n)\in\{(0,5),(0,6),(1,2),(2,0)\}$ are considered at the end of the paper.
It is also important to note that even in the case of five parabolic points on the Riemann sphere, a complete description of singular fibers is not known, see for example \cite[Discussion page 14]{S18}. In ongoing work, we shall use the tools developed in this paper to describe all the singular fibers of the Hitchin fibration when $(g_X,n)\in\{(0,5), (1,2)\}$. Here, the role of elementary transformations and the variation of the Hitchin fibration with respect to the two ramified elliptic cover are crucial.  As a last piece of motivation, the description of singular fibers includes  the study of the nilpotent cone and  the locus of fixed points of the $\C^*$-action given by multiplication on Higgs field. In particular, these techniques might also be useful to investigate the foliation conjecture \cite[Question 7.4]{S08}. Indeed, it follows from the work of C. T. Simpson that there is a decomposition in the moduli space of logarithmic connections obtained by looking at the limit in the moduli space of Higgs bundles, the limit point is a fixed point with respect to the $\C^*$-action and the foliation conjecture predicts that this decomposition forms a regular foliation.

Let $\omega_X$ denote the canonical sheaf of $X$. A rank two Higgs bundle on $(X, D)$ consists of a rank two holomorphic vector bundle $E$ on $X$,    endowed with  a  homomorphism $\theta: E \to E\otimes \omega_X(D)$, which has nilpotent residual matrix $Res (\theta; t_i)$ at each parabolic point. By $Res (\theta; t_i)$ we mean the linear endomorphism of the fiber $E_{t_i}$ defined by taking the residues at $t_i$ of local 1-forms defining $\theta$.
For any integer $d$, there is a moduli space $\mathcal H(X, D, d)$ parametrizing  Higgs bundles over $(X, D)$, with fixed degree $\deg E = d$. Actually, we consider triples $(E, { l}, \theta)$ where ${ l}$ intends to be the parabolic direction, described just below. See Section~\ref{PHB} for the detailed definition. The definition of these moduli spaces also depends on a choice of a weight vector (a sequence of reals numbers $0\le \mu_i \le 1$, $i=1,\dots, n$) which gives a notion of slope-stability.
The moduli space  $\mathcal H(X, D, d)$ is a normal variety of dimension $2(n-3+4g_X)$, if it is nonempty. For instance see \cite[Theorem 5.2]{Yo95} and our  Remark~\ref{rmk:open}.

Closely related, 
is the moduli space of parabolic vector bundles $\Bun(X,D, d)$, which is a normal variety of dimension $n-3+4g_X$, if the stable locus is nonempty. An element of it consists of a pair $(E, { l})$, where the additional data ${ l}$ refers to a one dimensional subspace on the fiber of $E$ over each parabolic point. There is a rational map $\mathcal H(X, D, d)\dashrightarrow \Bun(X, D, d)$ which associates to  a Higgs  field its underlying parabolic vector bundle, where the parabolic direction over $t_i$ coincides with the kernel of  $Res (\theta;t_i)$. It turns out that the fiber of this map, at a generic point, identifies with the cotangent space $T^*_{(E,{ l}, \theta)}  \Bun(X, D, d)$. In addition,  $\mathcal H(X, D, d)$ admits a symplectic structure which coincides with that of $T^*\Bun(X, D, d)$, given by the Liouville form.

An important player on this subject is the Hitchin map $H_X$,  it  associates to a Higgs bundle its characteristic polynomial.   This map is known to be an algebraically integrable system, which turns out to say that  the fiber over a generic point is Lagrangian and isomorphic to an Abelian variety. This last consists of the Picard  variety 
\[
H_X^{-1}(s)\simeq {\rm Pic}^{\frak n}(X_s) 
\]
(of the corresponding spectral curve $X_s$) which parametrizes line bundles of degree $\frak{n}=d+n+2(g_X-1)$ on $X_s$, by the  BNR correspondence. This correspondence is reviewed in Section~\ref{section:basic}, as well as the notion of spectral curves.

The main goal of this paper is to study the behavior of the Hitchin fibration under a degree two ramified covering $\pi: Y \to X$, between compact Riemann surfaces.  
We can pull back parabolic Higgs bundles from $X$ to $Y$.  The lifted Higgs field, after pull back, tends to have apparent singularities, as well as the corresponding spectral curve has many singular points. In order to eliminate them, we perform a suitable elementary transformation, and this process gives a correspondence $\Phi$ between moduli spaces, which we now describe.   

Assume that $\pi$ is branched over the points $x_1,\dots, x_n$, $n\ge 1$,  and let $B= \sum_{i=1}^n x_i$ be  the reduced divisor on $X$ defined by them.   Let $T= \sum_{i=1}^k t_i$, $k\ge 0$, be a reduced divisor on $X$ ($T=0$ when $k=0$) formed by points outside the support of  $B$.  We consider a map between moduli spaces
\[
\Phi: \mathcal H(X,B+T, 0) \dashrightarrow \mathcal H(Y, \pi^*T, -n)
\]
which is defined by doing the pull back $\pi^*$ followed by an elementary transformation $\elem_R$ over $R$, where $R = \sum_{i=1}^n y_i$ denotes the divisor on $Y$ defined by the set of ramification points, with $\pi(y_i)=x_i$. (see Section~\ref{elem_D} for the detailed definition). This map preserves the Hitchin fibration and we denote by 
\[
\Phi_s: H_X^{-1}(s) \to H_Y^{-1}(r)
\]
its restriction to a given fiber. Using the BNR correspondence, we obtain a map between the  Abelian varieties that parametrize line bundles on the spectral curves.  Hence, in order to describe the map $\Phi_s$ we are led to investigate the variation of the BNR correspondence with respect to elementary transformations, see Proposition~\ref{twistM}.     

Assuming $2g_X+n-4\ge 0$, with strict inequality if  $k=0$, we assure that the generic spectral curve in both moduli spaces is smooth and irreducible (see Proposition~\ref{prop:generalspeccurve}). We show  that there is an \'etale morphism $\xi_s: Y_r \to X_s$  of degree two between spectral curves (Proposition~\ref{prop:etale}), and the modular map $\Phi_s$ is determined by $\xi_s^*$, up to a translation. This is the main result of this paper,  Theorem~\ref{thm:main} in the main text.

\begin{thm}\label{thmintro:main}
Assume that  $2g_X+n-4\ge 0$, with strict inequality if  $k=0$. Then $\Phi$ is a rational map of degree two onto its image, and  it preserves the Hitchin fibrations.  Moreover, its restriction to a general Hitchin fiber is the map
\[
\Phi_s: {\rm Pic}^{\frak n}(X_s) \longrightarrow {\rm Pic}^{\tilde{\frak n}}({Y}_r)
\]
which sends a line bundle $M$ in ${\rm Pic}^{\frak n}(X_s)$ to the line bundle $\xi_s^*(M)({q}_r^*(-R))$ in ${\rm Pic}^{\tilde{\frak n}}({Y}_r)$.
\end{thm}

We finish the paper by applying Theorem~\ref{thmintro:main} to a degree two morphism $\pi: Y \to \P^1$ where $g_Y\in\{1, 2\}$. Let us briefly describe each case.

In the former case $g_Y=1$, we consider the moduli space $\mathcal H(\P^1, \Lambda)$ of $SL_2$-Higgs bundles  on $(\P^1, \Lambda)$, where $\Lambda$ is given by branch points $\{0, 1, \lambda, \infty\}$ of $\pi$ with an extra point $t$. By $SL_2$ we mean that Higgs fields  have vanishing trace and vector bundles have trivial determinant line bundle. After multiplication by a suitable line bundle to redress the determinant, we get a  modular map $\Phi_0:\mathcal H(\P^1, \Lambda)\to \mathcal H(Y, D)$  to the moduli space $\mathcal H(Y, D)$ of $SL_2$-Higgs bundles over the elliptic curve $Y$, with simple poles on $D=t_1+t_2$, with $\pi(t_i)=t$. The spectral curve $X_s$, on the $\P^1$ side, has genus $2$; on the elliptic side, the spectral curve $Y_r$ has genus $3$.   

In the second case $g_Y = 2$, we consider the moduli space $\mathcal H(\P^1, B)$ of $SL_2$-Higgs bundles over $(\P^1, B)$, where $B$ is formed by branch points of $\pi$. After renormalization of the determinant line bundle, we get a modular $\Phi_0:\mathcal H(\P^1, B)\to \mathcal H(Y)$ map with values on the moduli space $\mathcal H(Y)$ of holomorphic Higgs bundles (without parabolic points) over $Y$. The spectral curve $X_s$ over $\P^1$ has genus $3$, and over $Y$, the spectral curve $Y_r$ has genus $5$. In this context, Theorem~\ref{thmintro:main} yields the following result, which corresponds to Corollaries~\ref{cor:P1C} and \ref{cor:g=2} in the main text.

\begin{cor}
Let the notation be as above. In both cases $g_Y = 1, 2$, the modular map $\Phi_0$ is dominant and its restriction to a smooth Hitchin fiber is the degree two map 
\begin{eqnarray*}
\Phi_{0,s}: {\rm Pic}^{\frak n}(X_s) &\to& {\rm Prym}(Y_r/Y)\\ 
             M &\mapsto& \xi_s^*(M) \otimes q_r^*( L_0(- R))
\end{eqnarray*} 
where $L_0$ is a square root of $\mathcal O_Y(R)$ and ${\rm Prym}(Y_r/Y)$ is the Prym variety of the covering $q_r: Y_r\to Y$.   
\end{cor}

\section{Basic tools}\label{section:basic}

\subsection{Notations and conventions}\label{sec:notaandconv}

Throughout this paper, we assume that $X$ is a compact Riemann surface of genus $g_X\ge 0$. As usual, we denote by $\mathcal O_X$ its structural sheaf and  by $\omega_X$ its canonical sheaf.  When $L$ is a line bundle on $X$ we let $L^{-1}$ be its dual. If $E$ is a vector bundle on $X$ we often write $\Gamma (E)$ for ${\rm H}^0(X,E)$.

\subsection{Spectral curve}\label{sec:speccurve}
Let $L$ be a line bundle on $X$.  We denote by $\bf P$ the ruled surface $\P(\mathcal O_X\oplus L^{-1})$ which is the projectivization of the rank two vector bundle $\mathcal O_X\oplus L^{-1}$. Let   ${q}:{\bf P}\to X$ be the natural projection and let $\mathcal O_{\bf P}(1)$ be the hyperplane bundle along the fibers. We note that $q_*(\mathcal O_{\bf P}) = \mathcal O_X$ and $q_*(\mathcal O_{\bf P}(1)) = \mathcal O_X\oplus L^{-1}$, see for example \cite[II  - Proposition 7.11]{Ha}.

The ruled surface ${\bf P}$ contains two disjoint sections $\P(\mathcal O_X)$ and $\P(L^{-1})$ corresponding to the embedding of $\mathcal O_X$ and $L^{-1}$, respectively.  
The embedding   of $\mathcal O_X \hookrightarrow q_*(\mathcal O_{\bf P}(1))=\mathcal O_X\oplus L^{-1} $ gives a section ${\bf w}$ of $\mathcal O_{\bf P}(1)$, via the adjoint formula (cf. \cite[II - Section 5]{Ha})
which has $\mathbb P(L^{-1})$ as its zero set. Similarly, since 
\[
q_*(q^*L\otimes \mathcal O_{\bf P}(1))\simeq  L\otimes(\mathcal O_X\oplus L^{-1})=L\oplus \mathcal O_X  
\]
by the projection formula, we obtain a section ${\bf z}$ of $q^*L\otimes \mathcal O_{\bf P}(1)$ which has $\mathbb P(\mathcal O_X)$ as its zero set.

Given $s=(s_1,s_2)\in\Gamma(L)\oplus\Gamma(L^2)$, we define the {\it spectral curve} $X_s\subset {\bf P}$ as the zero locus  of the section 
\[
{\bf z}^2 +{q}^*(s_1) \cdot {\bf z}\cdot {\bf w} + { q}^*(s_2)\cdot {\bf w}^2 \in{\rm H}^0({\bf P}, {q}^*(L^2)\otimes \mathcal O_{\bf P}(2)). 
\]
It  comes with a degree two map  
\[
q_s:X_s\to X
\]
which is the restriction of $q$ to $X_s$. Notice that in a local open subset $U$ of $X$ where ${\bf P}|_U \simeq  U\times \mathbb P^1_{(w:z)}$, the spectral curve has equation
\begin{eqnarray*}\label{spectrallocal}
z^2+s_1zw+s_2w^2=0.
\end{eqnarray*}
Equivalently,   we can define a  structure of commutative ring  on  $\mathcal O_X\oplus L^{-1}$ induced by $s$:
\begin{eqnarray}\label{strucring}
(a_0,a_1)\cdot (b_0,b_1) := (a_0b_0-s_2a_1b_1,a_0b_1+a_1b_0-s_1a_1b_1).
\end{eqnarray}
This makes $\mathcal O_X\oplus L^{-1}$ an $\mathcal O_{X}$-algebra, which will be denoted by $\mathcal A_s$, and is locally given by 
\begin{eqnarray}\label{algebraOX}
\mathcal A_s(U) = \frac{\mathcal O_X(U)[z]}{(z^2 + s_1 z + s_2)}.
\end{eqnarray}
The spectral curve may be seen  as
\begin{eqnarray*}
X_s = {\rm Spec}~(\mathcal A_s).
\end{eqnarray*}
We often write simply 
\begin{eqnarray}\label{espcurves}
X_s = {\rm Spec}~(\mathcal O_X\oplus L^{-1})
\end{eqnarray}
where  $\mathcal O_X\oplus L^{-1}$ is considered as $\mathcal O_X$-algebra. Note  that 
\begin{eqnarray}\label{qsO}
(q_s)_*\mathcal O_{X_s} = \mathcal O_X\oplus L^{-1} .
\end{eqnarray}

The (arithmetic) genus of $X_s$ is given by 
\begin{eqnarray}\label{genusXs}
g_{X_s} := \dim {\rm H}^1(X_s,\mathcal O_{X_s}) = \deg L + 2(g_X-1)+1
\end{eqnarray}
see \cite[Remark 3.2]{BNR}.

\subsection{BNR correspondence}\label{correspondence}

Let us briefly recall the correspondence between  $\mathcal O_X$-homomorphisms $\theta: E \to E\otimes L$ having the same characteristic polynomial $P_s$ and  torsion free sheaves $M_{\theta}$ of rank one on $X_s$  (cf. \cite[Proposition 3.6]{BNR}). 

Let $E$ be a rank two holomorphic vector bundle on $X$ and let $\theta\in \Gamma(\End(E)\otimes L)$  be a homomorphism with $s_1 = -\tr(\theta)\in \Gamma(L)$ and $s_2 = \det (\theta)\in \Gamma(L^2)$. By Cayley-Hamilton theorem it satisfies the equation $P_s(\theta)=0$, where $P_s(z)=z^2+s_1z+s_2$ is the characteristic polynomial and 
\[
P_s(\theta) = \theta^2+s_1\cdot\theta +s_2\cdot I
\]
must be seen as a homomorphism $ E\to E\otimes L^2$. 

Let us assume that $X_s$ is integral, i.e. irreducible and reduced,  which turns out to say that $\mathcal A_s$ is a sheaf of  integral domains. We have a canonical isomorphism 
\[
\End(E)\otimes L \simeq \Hom_{\mathcal O_X}(L^{-1},\End(E))
\]
which associates $\theta\in \End(E)\otimes L$ to the $\mathcal O_X$-homomorphism
\begin{eqnarray*}
L^{-1} &\to& \End(E)\\
  z &\to& \theta_z
\end{eqnarray*}
where $\theta_z(e)=\theta(e)\cdot z$.
Thus the homomorphism $\theta$ induces a  structure of $\mathcal A_s$-module on $E$
\begin{eqnarray}\label{struc}
\Theta: \mathcal A_s \times E \to E
\end{eqnarray}
which is locally given by (see (\ref{algebraOX}))
\begin{eqnarray*}
\Theta_U: \mathcal A_s(U) \times E(U) &\to& E(U)\\
(a+bz, e) &\to& (aI+b\theta_z)\cdot e
\end{eqnarray*}
It defines a sheaf $M_{\theta}$ of $\mathcal A_s$-modules on $X_s$, which comes with a multiplication 
\begin{eqnarray}\label{structureofmodule}
(a+bz)\odot e := (aI+b\theta_z)\cdot e
\end{eqnarray}
for any $a+bz\in \mathcal A_s$, $e\in E$.  
Since $X_s$ is integral then $M_{\theta}$ is also torsion free.

Conversely,  if $M$ is a torsion free sheaf of rank one on $X_s$ then $E = (q_s)_*M$ is a locally free $\mathcal O_X$-module of rank two with an $\mathcal O_X$-linear map  
\[
\Theta: \mathcal A_s \times E \to E.
\]
The section ${\bf z}$ of $q_s^*L$ gives  a family of endomorphisms
\begin{eqnarray*}
\theta_U: E(U) &\to& E(U) \\
    s &\mapsto& \Theta(z,s)
\end{eqnarray*}
for each open set $U$ of $X$ such that ${\bf P}|_U \simeq  U\times \mathbb P^1_{(w:z)}$, which glue together to give a global homomorphism $\theta: E\to E\otimes L$ satisfying $\theta^2+s_1\theta+s_2I = 0 $. Since $X_s$ is integral then $P_s(z)=z^2+s_1z+s_2$ is irreducible over the function field of $X$. Thus $P_s$ is the characteristic polynomial of $\theta$.

Summarizing the discussion above, here is the version of BNR correspondence \cite[Proposition 3.6]{BNR} needed for our purposes:

\begin{prop}
Assume that the spectral curve $X_s$ is integral. Then there is a bijective correspondence between isomorphism classes of torsion free rank one sheaves on $X_s$ and isomorphism classes of pairs $(E, \theta)$ where $E$ has rank two and $\theta: E\to E\otimes L$ is an $\mathcal O_X$-homomorphism having characteristic polynomial $P_s$.  
\end{prop}

We conclude this section with a lemma which will be useful in the sequel.

\begin{lemma}\label{linebundle}
Let $M$ be the torsion free sheaf of rank one on $X_s$ associated to  $\theta: E \to E\otimes L$. Then $M$ is a line bundle on $X_s$ if and only if for every point $x\in X$ the $\mathbb C$-linear map $\theta_x: E_x \to E_x\otimes L_x\simeq E_x$ is not of the form $\lambda I$, where $\lambda\in \C$ and $I:E_x\to E_x$ is the identity. 
\end{lemma}

\proof
First assume  that $M$ is a line bundle on $X_s$, i.e. it has rank one as $\mathcal A_s$-module. Then for any $x\in X$  there exist an open subset $U\subset X$ containing $x$ and $v\in M(U)$ such that every $e\in M(U)$ writes as 
\[
e = (a+bz)\odot v, \quad a+bz\in \mathcal A_s(U).
\]
This means that any section $e\in E(U)$ can be written as 
\[
e = (a+bz)\odot v = av+b\theta(v), \quad a,b\in \mathcal O_X(U).
\]
In particular, $v$ and $\theta (v)$ are $\C$-linearly independent when considered as elements of $E_x$, which implies that $\theta_x \neq \lambda I$, for all $\lambda\in \C$. Reciprocally, assuming that $\theta_x \neq \lambda I$ for all $\lambda\in \C$,  there exist an open neighbourhood $U$ of $x$  and  $v\in E(U)$ such that $v$ and $\theta(v)$ are $\mathcal O_X(U)$-linearly independent.  This implies that $v\in M(U)$ generates $M(U)$ as $\mathcal A_s(U)$-module.
\endproof

\subsection{Elementary transformations}\label{elem_D}
We next introduce strongly parabolic endomorphisms and elementary transformations. The goal of this section is the behavior of $M_{\theta}$ (cf. Section \ref{correspondence}) with respect to elementary transformations. 

Fix  $t_1, \dots, t_n\in X$ distinct points and denote by $D=t_1+\cdots +t_n$ the  divisor determined by them. 

A \textit{quasiparabolic vector bundle} $ (E, { l})$, ${ l} = \{l_{i}\}$,  of rank two  on $\big(X, D\big)$ consists  of a holomorphic vector bundle $E$ of rank two on $X$ and for each  $i = 1,\dots,n$, a $1$-dimensional linear subspace $l_{i} \subset E_{t_{i}}$. We refer to the points $t_i$'s as parabolic points, and to the subspaces $l_{i} \subset E_{t_{i}}$ as the parabolic direction of $E$ at $t_i$.

There is a correspondence between quasiparabolic vector bundles, called elementary transformation,  which we now describe. We consider the following exact sequence of sheaves
$$
0\ \rightarrow\ E' \ \stackrel{\alpha}{\to} \  E\ \stackrel{\beta}{\rightarrow}\ \bigoplus_{i=1}^n E/l_i \ \rightarrow\ 0 \ 
$$
where $E/l_i$ intends to be a skyscraper sheaf determined by $E_{t_i}/l_i$,  i.e., for an open subset $U$ of $X$ we have  $ (E/l_i) (U) = E_{t_i}/l_i$ if $t_i\in U$ and $\{0\}$ otherwise. The map $\beta$ sends $s$ to $\oplus_{i=1}^n s(t_i)$. If $E$ is locally generated by $e_1,e_2$ as $\mathcal O_X$-module  near $t_i$ with $e_1(t_i)\in l_i$, then $E'$ is locally generated by $e_1,e_2'$, with $e_2' = xe_2$, where $x$ is a local coordinate. In particular $E'$ is locally free of rank two. We view $E'$ as a quasiparabolic vector bundle $(E',{ l'})$ of rank two over $(X,D)$ putting $l_i' := \ker\alpha_{t_i}$. We call it the {\it elementary transformation} of $(E,{ l})$ over $D$:
\[
\elem_D(E,{ l}) := (E', { l'}).
\]
Notice that we have the following equality
\[
\det E' = \det E \otimes \mathcal O_{X}(-D).
\]

An endomorphism $f : E \to E$ is called \textit{parabolic} if $f(l_{i}) \subseteq l_i$ for every $i = 1,\dots,n$. We denote by $\End(E, { l})$ the sheaf of parabolic endomorphisms.  A parabolic endomorphism is called \textit{strongly parabolic} if $f(E_{t_i}) \subseteq l_i$ and $f(l_i)=0$ for every $i = 1,\dots,n$. The sheaf of strongly parabolic endomorphisms of $E$ will be denoted by $\SEnd(E, { l})$.

Given a  parabolic homomorphism  $\theta: E\to E\otimes L$, i.e. 
\[
\theta\in\Gamma(\End(E, { l})\otimes L)
\]
we can perform an elementary transformation $\elem_D$ on the pair $( E,\theta)$, centered in ${ l}$. For instance, since $\theta$ is parabolic with 
respect to  $ l$, then $\theta(E')\subset E'\otimes L$ and its restriction induces a homomorphism 
\[
\theta': E' \to E'\otimes L
\]
which is parabolic with respect to the direction ${ l'}$ of $E'$. If $e_1$, $e_2$ are local sections which generate $E$ near a parabolic point $t_i$ with $e_1(t_i)\in l_i$ then $\theta$ is given by
\begin{eqnarray*}
\theta=\left(
\begin{array}{ccc} 
a & b  \\
xc & d  \\
\end{array}
\right) .
\end{eqnarray*}
Since $E'$ is locally generated by $e_1,e_2'$, with $e_2' = xe_2$, then the restriction of $\theta$ to $E'$ corresponds to the matrix
\begin{eqnarray*}
\theta'=\left(
\begin{array}{ccc} 
a & xb  \\
c & d  \\
\end{array}
\right).
\end{eqnarray*}
Note that trace and determinant do not change after an elementary transformation. 

The  condition of been strongly parabolic, i.e. $\theta\in\Gamma(\SEnd(E, { l})\otimes L)$ yields 
\begin{displaymath}
\left\{ \begin{array}{ll}
\tr(\theta)\in\Gamma(L(-D))\\
\det(\theta)\in\Gamma(L^2(-D))
\end{array} \right.
\end{displaymath}
and in this case the spectral curve $X_s$, determined by 
\[
s = (-\tr \theta, \det \theta)\in \Gamma(L)\oplus\Gamma(L^2)
\] 
is ramified over $D$.

Before stating the next result, we introduce some notation.  For this, let $\theta: E\to E\otimes L$ be a parabolic homomorphism with respect to ${ l}$ and assume that the corresponding spectral curve $X_s$ is  integral, in particular $\det \theta$ is nonzero. We shall consider the curve  $X_{\theta}\subset \P E$ of eigendirections of $\theta$, defined as follows. If we take a local trivialization  $\P E|_U\simeq U\times \P^1$ of the ruled surface $\P E$
at a point $x\in X$ where $\det \theta_x \neq 0$, and if $\theta$  has local matrix
\begin{eqnarray*}
\theta=\left(
\begin{array}{ccc} 
a & b  \\
c & d  \\
\end{array}
\right) 
\end{eqnarray*}
then we consider the curve in $\P E|_U$ having equation
\[
bz^2+(a-d)z-c = 0
\]
where $(1:z)$ is a coordinate for $\P^1$. These local equations patch together to build a quasi-projective subvariety of $\P E$ and $X_{\theta}$ is defined as the closure of this subvariety. Notice that $X_{\theta}$ intersects the fiber $\P E_x$, at a general point $x\in X$, exactly at two eigendirections of the linear homomorphism $\theta_x: E_x \to E_x\otimes L_x\simeq E_x$.

\begin{center}
\begin{figure}[h]
\centering
\includegraphics[height=2.8in]{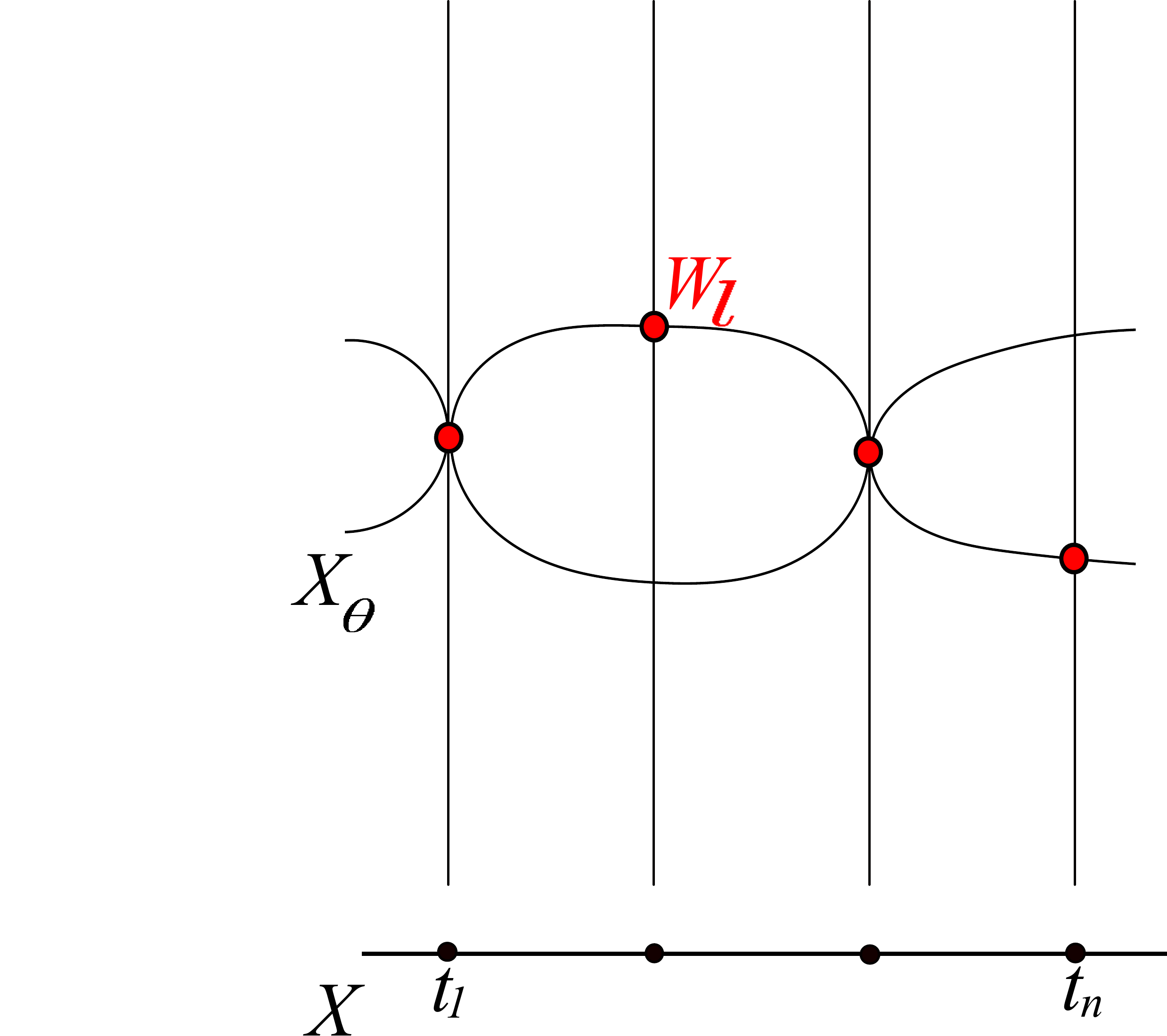}
\caption{Divisor $W_{ l}$ defined by parabolic directions.}
\label{divW}
\end{figure}
\end{center}

Let us denote by
\[
W_{ l} = l_1+\cdots + l_n
\] 
the {\it  divisor in $X_{\theta}$ defined by parabolic directions},  see Figure~\ref{divW}. There is a birational morphism 
\[
\eta: X_{\theta} \to X_s 
\]
which associates to each eigendirection the corresponding eigenvalue, and then it defines a divisor  
\[
W_{{ l}, s} = \eta(l_1)+\cdots +\eta(l_n)
\]
on $X_s$. If in addition $X_s$ is smooth, then $\eta$ is an isomorphism.  Note that when $\theta$ is strongly parabolic, and then $X_s$ is ramified over $D$, one has  $W_{{ l},s} = (q_s^*D)_{red}$. The subscript {\it red} denotes the effective reduced divisor defined by the support of $q_s^*D$.

The next result describes the behavior of $M_{\theta}$ (cf. Section~\ref{correspondence}) under elementary transformations.

\begin{prop}\label{twistM}
Assume that the spectral curve $X_s$ is integral. Let $M_{\theta}$ be a line bundle on $X_s$ corresponding to a  parabolic homomorphism  $\theta\in\Gamma(\mathcal{E}nd(E, { l})\otimes L)$. Let $M_{\theta'}$ be the torsion free sheaf corresponding to  $\theta'\in\Gamma(\mathcal{E}nd(E', { l}')\otimes L)$,  obtained by performing an elementary transformation over $D$. Then  
 \[
 M_{\theta'} = M_{\theta}\otimes \mathcal I
 \] 
where $\mathcal I$  is the ideal sheaf of the subscheme consisting of the support of $ q_s^*D - W_{{ l},s}$. 
\end{prop}

\proof
The homomorphism  $\theta: E\to E\otimes L$  induces a structure of $\mathcal O_{X_s}$-module 
\[
\Theta: \mathcal O_{X_s} \times E \to E
\]
on $E$ and  $M_{\theta}$ is just $E$ with this structure.  See Section \ref{correspondence}. Let $M_{\theta'}$ be the torsion free sheaf of rank one associated to $\theta'$. Since $\theta'$ is the restriction of $\theta$ to $E'$, then $M_{\theta'}$ is the $\mathcal O_{X_s}$-submodule of $M_{\theta}$ determined by the restriction $\Theta'$ of $\Theta$ to $E'$
\[
\Theta': \mathcal O_{X_s} \times E' \to E'.
\]

It is enough to analyse the behavior of $M_{\theta}$ at a point $t$ of the support of $D$, because outside it does not change. We may assume, without loss of generality, that $t=0$, $l_{0} = \begin{pmatrix}1\\ 0\end{pmatrix}$, and
\begin{eqnarray*}
\theta=\left(
\begin{array}{ccc} 
a & b  \\
xc & d  \\
\end{array}
\right)
\end{eqnarray*}
and the spectral curve $X_s$ has local equation 
\[
X_s = \{(x,z)\;\;: \;\; z^2-(a+d)z+(ad-xbc)=0\}.
\]
Since $M_{\theta}$ is a line bundle, then there is a generator $v$ for the $\mathcal O_{X_s}$-module $M_{\theta}$, i.e.,  any local section $e$ of $M_{\theta}$, can be written as    
\[
e  = (\alpha+\beta z)\odot v \in M_{\theta}
\]
with $\alpha+\beta z\in \mathcal O_{X_s}$   ($\alpha, \beta\in\mathcal O_X$, cf. (\ref{structureofmodule})), 
meaning that any local section $e$ of $E$   can be written as    
\[
e  = (\alpha I +\beta \theta) (v)  \in E.
\]
Sections of $E'$ are those of the form $e' = \alpha v+\beta \theta(v)\in E$, such that  
\[
e'_0  = (\alpha_0 I+\beta_0\theta_0)(v_0) 
\]
is proportional to $l_0= \begin{pmatrix}1\\ 0\end{pmatrix}$, and the subscript means evaluation at the point $x=0$. Since 
\begin{eqnarray*}
\alpha_0 I+\beta_0\theta_0=\left(
\begin{array}{ccc} 
\alpha_0+\beta_0a_0 & \beta_0b_0  \\
0 & \alpha_0+\beta_0d_0  \\
\end{array}
\right)
\end{eqnarray*}
one obtains that $\alpha_0+\beta_0d_0=0$, because $v_0$ is not proportional to $l_0$. This last implies that the regular function $\alpha+\beta z\in \mathcal O_{X_s}$ vanishes at the point $(0,d_0)\in X_s$, which corresponds to the second eigendirection of $\theta$ when $a_0\neq d_0$. When both eigendirections coincide, for instance when $X_s$ is singular or has a smooth ramification point, then $\alpha+\beta z$ vanishes at this point. 
 
We conclude that the sections of $M_{\theta'}$ are those sections of $M_{\theta}$ which vanish  at the support of the divisor $ q_s^*D - W_{{ l}, s} $. This gives $M_{\theta'}=M_{\theta}\otimes \mathcal I$ and this concludes the proof of the proposition. 
\endproof

\begin{cor}
With the same notation of Proposition~\ref{twistM}, if $X_s$ is smooth and irreducible,  and $\theta$ is strongly parabolic then 
\[
 M_{\theta'} = M_{\theta}( - (q_s^*D)_{red} ).
\]
\end{cor}

\proof
When $\theta$ is strongly parabolic,  $X_s$ is ramified over every point in the support of $D$, then $q_s^*D - W_{{ l},s} = (q_s^*D)_{red}$. 
In addition, since $X_s$ is smooth then $\mathcal I = \mathcal O_{X_s}(-(q_s^*D)_{red})$.
\endproof

\subsection{Pullback and twist}\label{pbtw}

Let $Y$ be another compact Riemann surface and assume that  $\pi: Y \to X$ is a finite morphism. 
Given $s=(s_1,s_2)\in\Gamma(L)\oplus\Gamma(L^2)$, let $Y_r^*$ be the spectral curve associated to $r=\pi^*s\in\Gamma(\pi^*L)\oplus\Gamma(\pi^*L^2)$.  The superscript $*$  in $Y_r^*$  will be justified in Section~\ref{sec:degtwo}, where a normalization $Y_r$ of $Y_r^*$ is constructed for our purposes.      

We note that  
\[
Y_r^*=Y\times_X X_s.
\]
Indeed, since $\pi$ is finite and thus affine,  we can assume $\pi:{\rm Spec}~B\to {\rm Spec}~A$, $X_s= {\rm Spec}~A_s$ and $Y_r^* = {\rm Spec}~B_{\pi^*s}$ with 
\[
A_s= \frac{A[z]}{(z^2+s_1z+s_2)}\;\;\text{and}\;\; B_{\pi^*s}=\frac{B[z]}{(z^2+\pi^*s_1z+\pi^*s_2)};
\]
and then we see that  $B_{\pi^*s} = B\otimes_A A_s$.   In terms of sheaves of algebras, 
we have 
\[
\mathcal B_{\pi^*s} = \mathcal O_Y \otimes_{\mathcal \pi^{-1}(\mathcal O_X)} \pi^{-1}(\mathcal A_s)
\]
where $\mathcal B_{\pi^*s}$ denotes the structural sheaf of $Y_r^*$.

In particular, there is a lifting $\pi_s:Y_r^*\to X_s$ making the following diagram commute
\begin{eqnarray}\label{diagpis}
\xymatrix { 
Y_r^* \ar@{->}[d] \ar@{->}[r]^{\pi_s}  &  X_s \ar@{->}[d] \\
Y \ar@{->}[r]^{\pi}    &       X 
}
\end{eqnarray}
and which is locally given by  the natural homomorphism $A_s\to A_s\otimes_A B=B_{\pi^*s}$. The next result describes the behavior of $M_{\theta}$ (cf. Section~\ref{correspondence})  under pullback.

\begin{prop}\label{prop:pullback}
Assume that both spectral curves  $X_s$ and $Y_r^*$ are integral.  Let $\theta: E\to E\otimes L$ be a homomorphism with characteristic coefficient $s$ and  let $M$ be the torsion free sheaf of rank one on $X_s$ associated to $\theta$. Then $\pi_s^*M$ is the torsion free sheaf of rank one associated to $\pi^*\theta$. 
\end{prop}

\proof
The pullback of $M$ by $\pi_s: Y_r^*\to X_s$ is given by 
\[
\pi_s^*M= \pi_s^{-1}(M)\otimes_{ \pi_s^{-1}(\mathcal O_{X_s})} \mathcal O_{Y_r^*}.
\]
If $U\subset Y$ and $V\subset X$ are local affine open sets where 
\[
\pi:U={\rm Spec}~B\to V={\rm Spec}~A
\] 
let us denote by $\tilde{U}\subset Y_r^*$ and $\tilde{V}\subset X_s$ their respective inverse images.    We can write 
\[
(\pi_s^*M)(\tilde{U}) = M(\tilde{V})\otimes_{A_s} B_{\pi^*s}. 
\]
Since the  $A_s$-module structure on $M(\tilde{V})$ is induced by $\theta$, then the $ B_{\pi^*s}$-module structure on $(\pi_s^*M)(\tilde{U})$ is induced by $\pi^*\theta$.   Which is equivalent to give a $ B_{\pi^*s}$-module structure on 
\[
\pi^*(E)(U) = E(V)\otimes_A B 
\]
induced by $\pi^*\theta$. It turns out to say that  $\pi_s^*M$ coincides with $\pi^*(E)$ with $ B_{\pi^*s}$-module structure induced by $\pi^*\theta$. This shows that  $\pi_s^*M$  is the torsion free sheaf of rank one associated to $\pi^*\theta$. 
\endproof

We finish this section by analyzing  the behavior of $M$ under twist.  

\begin{prop}\label{prop:twist}
Assume that $X_s$ is integral. Let $\theta: E\to E\otimes L$ be a homomorphism with characteristic coefficient $s$ and  let $M$ be the torsion free sheaf of rank one on $X_s$ associated to $\theta$. Let $L_0$ be a line bundle on $X$ and let $\theta_{L_0}: E\otimes L_0\to (E\otimes L_0)\otimes L$ be the homomorphism induced by $\theta$. Then the  torsion free sheaf of rank one on $X_s$ associated to $\theta_{L_0}$ is $M\otimes q_s^*L_0$. 
\end{prop}

\proof
Let $M_{L_0}$ be the torsion free sheaf of rank one on $X_s$ associated to $\theta_{L_0}$.  It follows from projection formula that 
\[
(q_s)_*(M\otimes q_s^*L_0) = ((q_s)_*(M))\otimes L_0 = E\otimes L_0
\]
and then $M\otimes q_s^*L_0$ and $M_{L_0}$ have the same pushforward via $q_s$. The reader can  check  that the structure of $\mathcal A_s$-module on $E\otimes L_0$ induced by $\theta_{L_0}$ coincides with that of  $M\otimes q_s^*L_0$.
\endproof

\section{Parabolic Higgs bundles and Hitchin fibration}

\subsection{Parabolic vector bundles}\label{pvb}

Let $(E, { l})$ be a quasiparabolic vector bundle. Fix a weight  vector $\mu = (\mu_{1}, \dots, \mu_{n})$ of real numbers $0 \leq \mu_{i} \leq 1$.
The  \textit{parabolic slope}  of $(E, { l})$ with respect to $\mu$ is 
$$
{\rm Slope^{\mu}}(E) = \frac{\deg E + \sum_{i=1}^{n}\mu_{i}}{2}
$$ 
where $\deg E = \deg (\det E)$.  Let $N \subset E$ be a line subbundle. For each  $i = 1,\dots,n$, set 
$${\rm Slope^{\mu}}_i(N,E)\ \ = \  \left\{ 
\begin{aligned}
& \mu_i \ & \text{ if } N_{t_{i}} = L_{i},\\
&0 \ &  \text{ if } N_{t_{i}} \neq L_{i}.
\end{aligned}
\right.$$
The  \textit{parabolic slope}  of $N \subset E$ with respect to $\mu$ is 
$$
{\rm Slope^{\mu}}(N,E) =  \deg N+\sum_{i=1}^{n} {\rm Slope^{\mu}}_i(N,E).
$$

A quasiparabolic vector bundle $(E,{ l})$  is $\mu$-\textit{semistable} (respectively $\mu$-\textit{stable}) if for every  line subbundle $N \subset E$ we have 
\[
{\rm Slope^{\mu}}(N,E) \leq  {\rm Slope^{\mu}}(E)
\] 
(respectively ${\rm Slope^{\mu}}(N,E) < {\rm Slope^{\mu}}(E)$). A \textit{parabolic vector bundle} is a quasiparabolic vector bundle together with a  weight vector $\mu$. We say that a parabolic vector bundle is {\it semistable} if the corresponding quasiparabolic vector bundle is $\mu$-\textit{semistable}.

For each  $d\in\mathbb Z$ and a weight vector $\mu$, there is a {\it moduli space} $\Bun_{\mu}(X,D,d)$, 
parametrizing rank two parabolic vector bundles $(E, { l})$ on $\big(X, D\big)$, with $\deg E = d$,  which are semistable.  If $N$ is a  line subbundle with $\deg N=d$, we denote by $\Bun_{\mu}(X,D,N)$ the subvariety of $\Bun_{\mu}(X, D, d)$ given by those parabolic vector bundles with $\det E = N$.

The moduli space $\Bun_{\mu}(X, D, d)$ has a structure of  a normal projective variety of dimension $n-3+4g_X$, if the stable locus  is not empty, see  \cite{MS80,Bho96}. By twisting vector bundles with a fixed line bundle $L_0$, we see that $\Bun_{\mu}(X, D, N) \cong \Bun_{\mu}(X, D, N\otimes L_0^2)$.

Elementary transformations give  correspondences between moduli spaces.  The stability condition is preserved after an appropriate 
modification of weights,  if $(E,{ l})$ is $\mu$-semistable then $\elem_D(E,{ l})$ is $\mu'$-semistable with
\[
\mu' = (1-\mu_1,\dots, 1-\mu_n).
\]
In particular, this gives  a map
$$
\elem_D\colon \Bun_{\mu}(X, D, N)   \to \Bun_{\mu'}(X, D, N\otimes \mathcal O_X(-D)).
$$

\subsection{Parabolic Higgs bundles}\label{PHB}

 A {\it parabolic Higgs bundle} is a triple $(E,{ l}, \theta)$ where $(E, { l})$ is a quasiparabolic vector bundle over $(X,D)$ and 
 $$
\theta \in \Gamma(\SEnd( E, { l})\otimes \omega_{X}(D))
$$
is a strongly parabolic homomorphism in the sense of Section~\ref{elem_D}. We say that $\theta$ is a {\it parabolic Higgs field}.  We note that been strongly parabolic means that for each point $t\in X$ lying in the support of $D$, the endomorphism 
$$
Res (\theta; t) \in {\rm End}(E_{t})
$$   
is nilpotent with respect to the parabolic direction of $E_{t}$. The nilpotency condition means that if $p \subset E_{t}$ is the parabolic direction, then 
\[
Res (\theta;t)(p) = 0 \quad \text{and} \quad Res (\theta; t)(E_{t})\subset p.
\] 
A line subbundle $N\subset E$ is called {\it invariant} under $\theta: E\to E\otimes\omega_X(D)$ if 
\[
\theta(N)\subset N\otimes \omega_X(D). 
\]
We say that $\theta$ is {\it irreducible} if it does not admit invariant line subbundle.

Fix a weight vector $\mu\in [0,1]^n$.  A parabolic Higgs bundle $(E, { l},\theta)$ on $\big(X, D\big)$ is called $\mu$-\textit{semistable} (respectively $\mu$-\textit{stable}) if for every  line subbundle $N \subset E$ invariant under $\theta$, we have 
\[
{\rm Slope^{\mu}}(N,E) \leq  {\rm Slope^{\mu}}(E)
\]
(respectively ${\rm Slope^{\mu}}(N,E) <  {\rm Slope^{\mu}}(E)$). We say that $(E, { l},\theta)$ is $\mu$-\textit{unstable} if it is not  $\mu$-\textit{semistable}. Note that an irreducible  Higgs field is stable for any choice of weight vector. 

Given $d\in \mathbb Z$, let ${\mathcal H}_{\mu}(X,D,d)$ be the moduli space of $\mu$-semistable parabolic Higgs bundles  $(E, { l},\theta)$ on $\big(X, D\big)$ with $\deg E=d$. 
If  $N$ is a line bundle of degree $d$, we denote by $\mathcal{H}_{\mu}(X,D,N)$ the subvariety of  $\mathcal{H}_{\mu}(X,D,d)$ given by those parabolic Higgs bundles  with fixed determinant $\det E=N$.

The moduli space $\mathcal H_{\mu}(X, D, d)$ is a quasiprojective variety of dimension $2(n-3+4g_X)$, see \cite{Yo93, Yo95, BY}. 

Similar to parabolic vector bundles, given  $(E, { l}, \theta)$ we can perform an elementary transformation $\elem_D$ on it, see Section \ref{elem_D}, and  this provides a  correspondence between moduli spaces of Higgs bundles
\[
\elem_D :  \mathcal{H}_{\mu}(X,D,N) \to \mathcal{H}_{\mu'}(X,D,N\otimes\mathcal O_X(-D)).
\]

\subsection{Hitchin map and its fibers} \label{section Hitchin}
Let $(E, { l},\theta)\in {\mathcal H}_{\mu}(X,D,d)$ be a parabolic Higgs bundle.
Since $Res (\theta;t_i)$ is nilpotent for every parabolic point $t_i\in X$, one obtains 
\begin{displaymath}
\left\{ \begin{array}{ll}
\tr(\theta)\in\Gamma(\omega_{X})\\
\det(\theta)\in\Gamma(\omega_{X}^{\otimes 2}(D))\subset \Gamma(\omega_{X}^{\otimes 2}(2D)).
\end{array} \right.
\end{displaymath}

In order to simplify notation we shall write $\Sigma_D$ for the vector space 
\begin{eqnarray}\label{defSigma}
\Sigma_D = \Gamma(\omega_{X})\oplus\Gamma(\omega_{X}^{\otimes 2}(D)). 
\end{eqnarray}
The \textit{Hitchin map} is defined as
$$
\begin{array}{cccc}
H_X: &{\mathcal H}_{\mu}(X,D,d)& \longrightarrow & \Sigma_D\\
      & (E,{ l},\theta) & \longmapsto & (-\tr(\theta), \det(\theta)).
\end{array}
$$

Let $s\in \Sigma_D$. Using the correspondence of Section \ref{correspondence} we get an identification between the Jacobian variety of $X_s$ and the fiber of $H_X$ over $s$, when $X_s$ is integral. Then, we are led to investigate the existence of integral spectral curves. To do this, we need a couple of results.

\begin{prop}\label{non-integral}
Let $(E, { l},\theta)$ be a parabolic Higgs bundle on $\big(X, D\big)$ and let $s=(-\tr\theta, \det \theta)$.  If the spectral curve $X_s$ is non-integral  then $\det \theta \in \Gamma(\omega_X^{\otimes 2})$. 
\end{prop}

\proof
If $X_s$ is non-integral then  the characteristic polynomial $P_s(z) = z^2  - \tr\theta z +\det \theta$ is reducible over the function field of $X$. Hence, near  a parabolic point $t_i$ in the support of $D$ we can write 
\[
P_s(z) = (z-a)(z-b)
\]
with $a, b\in \mathcal O_X(U)$. The nilpotence condition on $\theta$  implies that $a+b$ and $ab$ vanish at $t_i$. In particular, $a$ and $b$ vanish at $t_i$ and then $\det \theta\in\Gamma(\omega_X^{\otimes 2}(2D))$ has a zero of order at least two at $t_i$, $i=1,\dots, n$, which means that $\det \theta \in \Gamma(\omega_X^{\otimes 2})$. 
\endproof

\begin{prop}\label{prop:invar}
Let $(E, { l},\theta)$ be a parabolic Higgs bundle on $\big(X, D\big)$. If $\theta$ has an invariant line bundle $N\subset E$,   then $X_s$ is non-integral.
\end{prop}

\proof
We will show that the characteristic polynomial $P_s(z) = z^2  - \tr\theta z +\det \theta$ is reducible over the function field of $X$. Taking a Zariski open subset $U$ of $X$ and $e\in N(U)$ there exists $\lambda\in \omega_X(D)(U)$ such that 
\[
\theta(e) = \lambda \cdot e.
\]
In particular, the rational function $z-\lambda$ on $X$  divides  $P_s(z)$.
\endproof

Let us  assume that  $X_s$ is integral   
and let $M$ be a line bundle on $X_s$. We can compute the degree of $E=(q_s)_*M$ using (\ref{qsO}) and the following identity 
\begin{eqnarray}\label{detE}
{\rm det} (E) \simeq {\rm det}((q_s)_*\mathcal O_{X_s} )\otimes {\rm Nm}(M)
\end{eqnarray}
where ${\rm Nm}(M)$ is the norm map  (see \cite[Cap. IV Ex. 2.6]{Ha}). 
The result is 
$$
\deg M = \deg E+\deg L 
$$
where $L = \omega_X(D)$, and this yields 
$$
\deg M = \deg E + \deg D +2g_X - 2.
$$
This computation can be extended to torsion free rank one sheaves on $X_s$. Therefore, given $d\in \mathbb Z$ we are led to define 
\[
{\frak  n} = {\frak  n}(d) := d+ \deg D +2g_X-2. 
\]

Let us denote by $\overline{\rm Pic}^{\frak n}(X_s)$ the variety parametrizing isomorphism classes of   torsion free sheaves of rank one on $X_s$ of degree ${\frak  n}$ and let ${\rm Pic}^{\frak n}(X_s)$ denote its corresponding subset formed by line bundles.  For an irreducible curve having only planar singularities, $\overline{\rm Pic}^{\frak n}(X_s)$ is an irreducible variety that contains ${\rm Pic}^{\frak n}(X_s)$ as a dense open subset (c.f. \cite{AIK77}). If $X_s$ has a non-planar singularity then  $\overline{\rm Pic}^{\frak n}(X_s)$ has at least two irreducible components (c.f. \cite{KK81}).

As we have seen in Section \ref{correspondence},  there is a bijective correspondence between $\overline{\rm Pic}^{\frak  n}(X_s)$ and isomorphism classes of pairs $(E, \theta)$ where $E$ is a vector bundle of rank two and degree $d$, and $\theta: E\to E\otimes \omega_X(D)$ is a homomorphism with $\tr(\theta)=-s_1$ and $\det (\theta) = s_2$. In particular we have the following result (see also \cite{LM}).

\begin{prop}\label{HitFibre}
Assuming that $s = (s_1, s_2)\in \Sigma_D$, with
\[
s_2\in\Gamma(\omega_X^2(D))\setminus \Gamma(\omega_X^2(D-t_i))\;,\; \forall i=1, \dots, n
\]
then $X_s$ is integral and there is a bijective correspondence 
 \[
 \overline{\rm Pic}^{\frak n}(X_s) \longleftrightarrow H_X^{-1}(s)\;.
 \]
 Moreover, any Higgs bundle in $H_X^{-1}(s)$ has  irreducible Higgs field.
\end{prop}

\proof
Our hypothesis on $s_2$ ensures, in particular,  that $s_2\notin \Gamma(\omega_X^{\otimes 2})$, then $X_s$ is integral by Proposition \ref{non-integral}. 

Let $L=\omega_X(D)$. For each  $M\in \overline{\rm Pic}^{\frak n}(X_s)$ we can associate a pair $(E,\theta)$, where $\theta: E\to E\otimes L$ is a homomorphism  with characteristic polynomial $s$, see Section \ref{correspondence}.  Since $s_2$ lies in the subspace $\Gamma(L^2(-D))$ of $\Gamma(L^2)$ then any residual matrix $Res (\theta;t_i)$  has vanishing determinant. This and  $\tr (\theta)\in \Gamma(L(-D))$ yield  $\theta$  strongly parabolic.  Note that $Res (\theta;t_i)$  is non-null because $s_2\notin \Gamma(L^2(-D-t_i))$. Thus the parabolic direction $l_i\subset E_{t_i}$ is defined as the kernel of $Res (\theta;t_i)$. Now we  show that  $(E, { l}, \theta)$ is $\mu$-semistable. Actually, we can show that $\theta$ is irreducible and this also completes the last assertion of the statement of the proposition. Indeed,  by Proposition~\ref{prop:invar},  if there is a line subbundle $N\subset E$ invariant by $\theta$, then  $\det(\theta)=s_2\in\Gamma(\omega_X^{\otimes 2})$ and this gives  a  contradiction.

Conversely, for each pair $(E,\theta)$ with characteristic coefficient $s$,  BNR correspondence gives an element $M\in \overline{\rm Pic}^{\frak n}(X_s)$. This finishes the proof of the proposition. 
\endproof

\begin{remark}\label{rmk:open}\rm
Let $\mathcal U = \Gamma(\omega_X^2(D))\setminus \cup_{i=1}^n\Gamma(\omega_X^2(D-t_i))$ be the open set given by Proposition \ref{HitFibre}.  If $2g_X+\deg D- 4 \ge 0$ then $\mathcal U$ is nonempty. To see this, note that this inequality yields $\deg \omega_X^2(D)\ge 2g_X$ and then $\omega_X^2(D)$ has no base points. In particular, by Proposition~\ref{HitFibre},  if $2g_X+\deg D- 4 \ge 0$ then there exist irreducible parabolic  Higgs Higgs fields $\theta: E\to E\otimes\omega_X(D)$, and therefore $\mathcal H_{\mu}(X, D, d)$ contains stable Higgs bundles for any choice of weight vector $\mu$.
\hfill$\lrcorner$  
\end{remark}

We finish this section by considering  traceless Higgs fields.  In view of Lemma~\ref{linebundle} and the next result, Higgs fields having non-vanishing residual part at every parabolic point play a key role.  We say that $\theta: E\to E\otimes \omega_X(D)$ is {\it nowhere-holomorphic in} $D$ if 
\[
Res (\theta;t_i)\neq 0
\] 
for every $t_i$ at the support of $D$.

\begin{prop}\label{prop:lineversus}
Let $s=(0,s_2)\in\Sigma_D$. Assume that $X_s$ is integral and smooth at every point over the complement of the support of $D$. Then there is a bijective correspondence
\[
{\rm Pic}^{\frak n}(X_s) \longleftrightarrow \left\{ (E, { l}, \theta)\in H_X^{-1}(s)\;:\; \theta \;\text{is nowhere-holomorphic in}\;D\right\}
\]
\end{prop}

\proof
By Lemma~\ref{linebundle} there is a bijective correspondence between ${\rm Pic}^{\frak n}(X_s)$ and the isomorphism classes of pairs $(E,\theta)$ where $E$ has degree $d$, $\theta: E\to E\otimes \omega_X(D)$ is a homomorphism with $\tr(\theta)=0$ and $\det (\theta) = s_2$, and $\theta_x$ is not of the form $\lambda I$ for every $x\in X$. 

Since $\tr\theta=0$, then $\theta_x=\lambda I$ if and only if $\theta_x=0$. But, $\theta_x$ is nonvanishing at every point over the complement of the support of $D$, because $X_s$ is smooth at that point. For a point $t_i$ of the support of $D$, we have 
\[
\theta_{t_i} = Res (\theta;t_i)
\]
and then $\theta_{t_i}\neq 0$, $\forall i=1,\dots, n$, if and only if $\theta$ is nowhere-holomorphic in $D$.  To conclude the proof of the proposition we note that, in this case,  the parabolic direction over $t_i$ is determined by the kernel of $Res (\theta;t_i)$. 
\endproof

\section{Spectral curves under degree two ramified coverings}\label{sec:degtwo}

Let $\pi: Y \to X$ be a finite morphism of degree two between irreducible smooth complex curves. Assume that $\{x_1,\dots, x_n\}$, $n\ge 1$, is the set of branch points and let $B= \sum_{i=1}^n x_i$ be  the reduced divisor on $X$ defined by them. Let $R = \sum_{i=1}^n y_i$ be the divisor on $Y$ formed by ramification points: $\pi(y_i)=x_i$. 

Let us consider a reduced divisor $T= \sum_{i=1}^k t_i$, $k\ge 0$,  on $X$  formed by points outside the support of  $B$. We put $T=0$ when $k=0$. Then we can write 
\[
\pi^*(T+B) = \pi^*(T)+2R. 
\]
Note that $\pi$ induces a linear map (see Section~\ref{section Hitchin} for the definition of $\Sigma_D$)
\[
\pi^*: \Sigma_{T+B} \to \Sigma_{\pi^*(T)},
\]
because  poles over branch points became regular after pullback. Indeed, if  $y$ and $x$ denote local coordinates of $Y$ and $X$, respectively, over a branch point,  then  $\pi$ is locally given by $\pi(y)=y^2$ and 
\[
\pi^*\left(\frac{dx^{\otimes 2}}{x}\right) = 4dy^{\otimes 2}.
\]

We want to pull back Higgs fields over  $(X, T+B)$ to Higgs fields over $(Y, \pi^*(T))$. Let us first investigate the relation between their spectral curves.   

\subsection{General spectral curve}\label{sec:gens}

In this section we start with an element $s\in \Sigma_{T+B}$, which gives a spectral curve $X_s$. Since $\pi^*(s)$ becomes regular over ramification points, we can see it as an element  $r=\pi^*(s)\in\Sigma_{\pi^*(T)}$. Let $Y_r$ be the corresponding spectral curve given by $r$.  

We set $L=\omega_X(T+B)$  and consider the respective spectral curve (cf.~(\ref{espcurves}))
\[
X_s = {\rm Spec}~(\mathcal O_X\oplus L^{-1}) 
\] 
with its degree two map $q_s:X_s \to X$. The divisor $R_{q_s}$ of ramification points of $q_s$ is given by zeros of $s_1^2- 4s_2$ and since $s$ lies in the subspace $\Sigma_{T+B}$, then $q_s$ is ramified over every point in the support of $T+B$.

\begin{center}
\begin{figure}[h]
\centering
\includegraphics[height=1.3in]{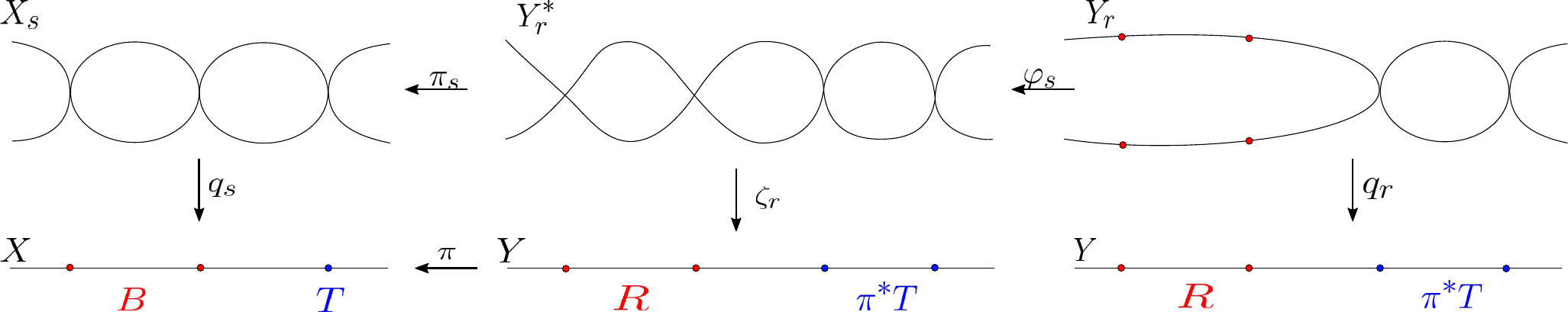}
\caption{Three spectral curves.}
\label{fig:especcurves}
\end{figure}
\end{center}

Note that 
\[
\pi^*L \simeq \omega_Y(R+\pi^*(T)).
\]
Following Section~\ref{pbtw}, the element 
\[
r=\pi^*s\in \Gamma(\pi^*L)\otimes\Gamma(\pi^*L^2)
\] 
gives the spectral curve 
\[
Y_r^* = {\rm Spec}~(\mathcal O_Y\oplus (\pi^*L)^{-1}) = Y\times_X X_s
\] 
which  is singular over a point lying in the support of  $R$. 
To avoid singular points of $Y_r^*$, it is convenient to switch $\pi^*L$. Then we  consider 
\[
N = \pi^*L(-R) = \omega_Y(\pi^*T)
\]
and the element 
\[
r=\pi^*s \in \Sigma_{\pi^*(T)}\subset \Gamma(N)\oplus \Gamma(N^2)
\]
gives the spectral curve 
\[ 
Y_r= {\rm Spec}(\mathcal O_Y\oplus N^{-1}). 
\] 
We summarize the discussion above in Figure~\ref{fig:especcurves}.

\begin{prop}\label{prop:generalspeccurve}
Assume that  $2g_X+n-4\ge 0$, with strict inequality if  $k=0$ ($n=\deg B$ and $k=\deg T$).  If $s$ is a general element of $\Sigma_{B+T}$, then $X_s$ and $Y_r$ are irreducible and smooth. 
\end{prop}

\proof
The hypothesis $2g_X+n-4\ge 0$ ensures that $X_s$ is integral for general $s\in \Sigma_{B+T}$,  see Remark~\ref{rmk:open}. 
The locus of sections $s\in \Sigma_{B+T}$, $s=(s_1, s_2)$,  such that $X_s$ is smooth is open and nonempty. Besides that if $s_2$ has only simple  zeros then the same holds for $\pi^*(s_2)\in\Sigma_{\pi^*(T)}$. In particular, the curve $Y_r$, $r=\pi^*(s)$,  is smooth for general $s$. Now we show that  $Y_r$ is integral. 

First let us assume  $k>0$ and  let $y$ be a point  with $\pi(y)=t_i$, where $t_i$ lies  in the support of $T$. If  $s_2\notin \Gamma(\omega_X^{\otimes 2}(B+T-t_i))$ then $\pi^*s_2\notin \Gamma(\omega_Y^{\otimes 2}(\pi^*T-y))$  because $\pi$ is \'etale at $y$. Therefore we can use Proposition \ref{HitFibre} again to conclude that $Y_r$ is integral.

Now  suppose  $k=0$. 
We claim that the image of  $\pi^*: \Sigma_{B+T}\to \Sigma_{\pi^*{(T)}}$ intersects the open subset of $\Sigma_{\pi^*{(T)}}$  formed by all sections for which the spectral curve is integral. To do this, we can assume $s_1=0$.  
A section $s_2\in \Gamma(\omega_Y^{\otimes 2})$ gives a reducible spectral curve if and only if it lies in the image of the map $\psi: \Gamma(\omega_Y)\to \Gamma(\omega_Y^{\otimes 2})$,  $a\mapsto a^2$. By Riemann-Roch theorem, the injective linear map $\pi^*: \Gamma(\omega_X^{\otimes 2}(B) )\to \Gamma(\omega_Y^{\otimes 2})$ has image of dimension $3(g_X-1)+n$, whereas the image of $\psi$ has dimension at most $g_Y$. Thus  if $3(g_X-1)+n>g_Y$, the claim above will follow. By Riemann-Hurwitz formula,   the previous inequality is equivalent to $2g_X+n-4>0$. This concludes the proof of the proposition.  
\endproof

\subsection{A map between spectral curves}

\begin{prop}\label{prop:etale}
 Let the notation be as above. The following statements hold: 
\begin{enumerate}
\item\label{existence} there is a morphism $\xi_s:Y_r\to X_s$ of degree two making the following diagram commute
 \[
\xymatrix { 
Y_{r} \ar@{->}[d]_{q_r} \ar@{->}[r]^{\xi_s}  &  X_s \ar@{->}[d]^{q_s} \\
Y \ar@{->}[r]^{\pi}    &       X \;.
}
\]
\item\label{normalization} if $X_s$ is  smooth then  $Y_r$  is the normalization of $Y_r^*$; 
\item\label{etale} $\xi_s$ is \'etale when $X_s$ and  $Y_r$ are smooth and irreducible. 
\end{enumerate}
\end{prop}

\proof

Let us prove (\ref{existence}).
Following Section~\ref{sec:gens}, setting $L=\omega_X(T+B)$ and $N = \omega_Y(\pi^*T)$ we have
\begin{displaymath}
\left\{ \begin{array}{ll}
X_s = {\rm Spec}~(\mathcal O_X\oplus L^{-1})\\
Y_r^* = {\rm Spec}~(\mathcal O_Y\oplus (\pi^*L)^{-1}) = Y\times_X X_s\\
Y_r = {\rm Spec}(\mathcal O_Y\oplus N^{-1}). 
\end{array} \right.
\end{displaymath}
Since  $(\pi^*L)^{-1}\subset N^{-1}$, this gives an inclusion  $i :\mathcal O_Y\oplus (\pi^*L)^{-1} \to \mathcal O_Y\oplus N^{-1}$ of $\mathcal O_Y$-modules. Note that  their structure of commutative ring induced by $\pi^*s$ (c.f.  (\ref{strucring})) are compatible and thus  $\mathcal O_Y\oplus (\pi^*L)^{-1}$ is  a subring of $\mathcal O_Y\oplus N^{-1}$.  In particular $i$ yields a  map $\varphi_s: Y_r \to Y_r^*$. Then  we set  $\xi_s = \pi_s\circ \varphi_s$, where $\pi_s: Y_r^*\to X_s$ is given by (\ref{diagpis}). Note that $(\pi^*L)^{-1}$ and $N^{-1}$ coincide when restricted to an open set $U$ away from the support of $R$, then $\varphi_s|_U$  is an isomorphism. Since $\pi_s:Y^*_r\to X_s$ has degree two and $\varphi_s: Y_r \to Y_r^*$ is birational, this concludes the proof of (\ref{existence}). 

In order to prove (\ref{normalization}), we note first that if  $X_s$  is smooth then  $s_2$ has only simple zeros. Consequently, the same is true for $\pi^*s_2$, seen  as a section of $\Gamma(N^2)$. This shows that  $Y_r$ is smooth and then $\varphi_s: Y_r \to Y_r^*$ is the normalization of $Y_r^*$. This proves (\ref{normalization}). 

Finally, when $X_s$ and $Y_r$ are smooth and irreducible,   we can use (\ref{genusXs}) and Riemann-Hurwitz formula to conclude (\ref{etale}).  
\endproof

In view  of  Proposition~\ref{prop:etale} we conclude this section with the following well-known result, see also \cite[Lemma p. 332]{Mu}.  Here, we include a proof for the convenience to the reader.

 \begin{prop}\label{prop:etalegeral}
 Let $\xi:Z\to W$ be an \'etale morphism of degree $2$ between compact Riemann surfaces. Then the induced morphism $ \xi^*: {\rm Pic}^0(W) \to {\rm Pic}^0(Z)$ has degree $2$ onto its image. 
 \end{prop}

\proof

There is a  correspondence between \'etale morphisms of degree $2$ and $2$-torsion elements in  ${\rm Pic}^0(W)$. 
We  recall briefly this correspondence. Let $L\in {\rm Pic}^0(W)$ be a $2$-torsion point, ${\bf P}=\P(\mathcal O_W\oplus L)$ and $q:{\bf P}\to W$ the natural projection. Roughly speaking, the curve $Z$ can be seen as the inverse image of the unit section of $L^2\simeq \mathcal O_W$ under the map $L\to L^2$ which sends $a$ to $a^2$.  The subvariety $\P(\mathcal O_W)\subset {\bf P}$ is the zero locus of a section 
$$
{\bf z}:\mathcal O_{\bf P} \to q^*(L)\otimes \mathcal O_{\bf P}(1)
$$ 
and $\P(L)$ is the zero locus of a section 
$$
{\bf w}:\mathcal O_{\bf P} \to \mathcal O_{\bf P}(1)
$$
and $\xi=q|_Z$. 
If we denote by $s$ a nonzero section of $L^2$, the curve $Z\subset {\bf P}$ corresponds to the zero locus  of the section 
\[
{\bf z}^2 - { q}^*(s)\cdot {\bf w}^2.
\]
Note that $s$ is not a power $a^2$ of a section $a$ of $L$, otherwise $L$ would be trivial.     

Now we can show that $\xi^*(L)= \mathcal O_Z$. To see this, observe that  ${\bf z}$ and ${\bf w}$ do not vanish in $Z$, then the restriction of $ \mathcal O_{\bf P}(1)$ and $q^*(L)\otimes \mathcal O_{\bf P}(1)$ to $Z$ are trivial.   It follows  that $$\xi^*(L)  = \left(q^*(L)\otimes \mathcal O_{\bf P}(1)\right)|_Z$$ is trivial. 

Assume that $M\in {\rm Pic}^0(W)$ is a  line bundle with $\xi^*(M) = \mathcal O_Z$.  Then projection formula yields 
$$
\mathcal O_W\oplus L \simeq \xi_*(\mathcal O_Z) = \xi_*((\pi^*M)\otimes \mathcal O_Z)\simeq M\oplus (M\otimes L).
$$
In particular, $M\oplus (M\otimes L)$ has a nonzero section. Since $M$ and $L$ have degree zero, one obtains that either $M$ is  $\mathcal O_W$  or $M\otimes L$ is $\mathcal O_W$. This concludes the proof of the proposition.  
\endproof

\section{A map between moduli spaces}

\subsection{Main result}\label{def phi} 
We keep the notation  of Section~\ref{sec:degtwo}, where $B$ and $R$ denote the divisors of branch points and ramification points of $\pi: Y\to X$, respectively. Also, let $T$ be an effective reduced divisor on $X$ whose support does not intersect the support of $B$.

In this section we will consider a rational map between moduli spaces of parabolic Higgs bundles 
\begin{eqnarray}\label{mapphimodulaire}
\Phi: \mathcal{H}_{\mu}(X,B+T,0) \dashrightarrow \mathcal{H}_{\mu'}(Y,\pi^*(T),-n).
\end{eqnarray}
 This map is defined by doing $\pi^*$ followed by an elementary transformation over $R$
\[
\xymatrix { {(\mathcal E,\theta)} \ar@{->}[r]_{\pi^*} \ar@/^0.5cm/[rr]^{\Phi}  &  (\pi^*\mathcal E, \pi^*\theta) \ar@{->}[r]_{\elem_R} & ((\pi^*\mathcal E)',( \pi^*\theta)')}
\]
where $\mathcal E = (E, { l})$. In the last step, parabolic directions over ramification points are forgotten, this justifies  the lack of ramification points in the divisor of parabolic points on the target in (\ref{mapphimodulaire}).  We describe below each step of this correspondence.

At first, $\Phi$ might send each semistable Higgs bundle to an unstable Higgs bundle. So, in view of Proposition~\ref{prop:generalspeccurve}, we assume that  $2g_X+n-4\ge 0$, with strict inequality if  $k=0$ ($n=\deg B$ and $k=\deg T$). Let us fix $s\in\Sigma_{B+T}$ such that $X_s$ and $Y_r$ are integral.

We can define a correspondence 
\[
(E,{ l},\theta) \mapsto \Phi (E,{ l},\theta)
\]
which associates to a Higgs bundle on $(X, T+B)$,  with characteristic coefficient $s$,  a Higgs bundle on $(Y, \pi^*(T))$ with characteristic coefficient $r=\pi^*s$, by doing the following series of transformations:
\begin{enumerate}
\item Let $(E,{ l},\theta)$  be a parabolic Higgs bundle, with $\deg E = 0$, and with characteristic coefficient $s$. 
\item Pulling back $(E,{ l},\theta)$ to $Y$ via $\pi$,  we obtain a parabolic Higgs bundle $(\pi^*(E,{ l}),\pi^*\theta)$ on $(Y, \pi^*(B+T))$ with characteristic coefficient $r\in \Sigma_{\pi^*(T)}$. 
\item After performing an elementary transformation over $R$ (cf. Section~\ref{elem_D}), the resulting parabolic Higgs bundle  
\[
\Phi(E,{ l}, \theta):= \big(((\pi^*E)', { l}'),(\pi^*\theta)'\big)
\]  
is regular over each ramification point, and $\deg (\pi^*E)' = -n$. Then we get a parabolic Higgs bundle over $(Y, \pi^*(T))$. Since $X_s$ and $Y_r$ are integral, $\theta$ and  $(\pi^*\theta)'$ are irreducible (cf. Proposition~\ref{prop:invar}), then they are both stable.  
\end{enumerate}

This correspondence gives a rational map between moduli spaces of parabolic Higgs bundles
\begin{eqnarray}\label{themap}
\Phi: \mathcal{H}_{\mu}(X,B+T,0) \dashrightarrow \mathcal{H}_{\mu'}(Y,\pi^*(T),-n)
\end{eqnarray}
which makes the following diagram commute 
 \begin{eqnarray}\label{diag}
\xymatrix { 
 \mathcal{H}_{\mu}(X,B+T, 0) \ar@{->}[d]_{H_X} \ar@{-->}[r]^{\Phi}  & \mathcal{H}_{\mu'}(Y,\pi^*(T), -n) \ar@{->}[d]^{H_Y} \\
\Sigma_{B+T} \ar@{->}[r]^{\pi^*}    &        \Sigma_{\pi^*(T)}
}
\end{eqnarray}

We want to investigate the behavior of $\Phi$ with respect to the Hitchin fibration. 
Let us consider  its restriction to a (general) Hitchin fiber
\[
\Phi_s : H_X^{-1}(s) \to H_Y^{-1}(r)
\]
with $r = \pi^*s$. When $X_s$ and $Y_r$ are irreducible and smooth, then 
\[
H_X^{-1}(s) \simeq {\rm Pic}^{\frak n}(X_s) \quad\text{and}\quad H_Y^{-1}(r) \simeq {\rm Pic}^{\tilde{\frak n}}(Y_r).
\]
where ${\frak  n} = k+n+2g_X-2$ and ${\tilde{\frak n}} = 2k-n+2g_Y-2$. This gives a map 
\[
\Phi_{s}: {\rm Pic}^{\frak n}(X_s) \to  {\rm Pic}^{\tilde{\frak n}}(Y_r).
\]

In the next result we show that this map is determined by $\xi_s^*$ (see Proposition~\ref{prop:etale}), up to a translation on the target.

\begin{thm}\label{thm:main}
Assume that  $2g_X+n-4\ge 0$, with strict inequality if  $k=0$. Then $\Phi$ is a rational map of degree two onto its image which preserves the Hitchin fibrations.  Moreover, its restriction to a  general Hitchin fiber is the map
\[
\Phi_s: {\rm Pic}^{\frak n}(X_s) \longrightarrow  {\rm Pic}^{\tilde{\frak n}}(Y_r)
\]
which sends a line bundle $M$ in ${\rm Pic}^{\frak n}(X_s)$ to the line bundle $\xi_s^*(M)(q_r^*(-R))$ in ${\rm Pic}^{\tilde{\frak n}}(Y_r)$, 
where $\xi_s$ is given by Proposition~\ref{prop:etale}.
\end{thm}

\proof
By Proposition~\ref{prop:generalspeccurve}, $X_s$ and $Y_r$ are irreducible and smooth for general $s$.   We want to show that $\Phi_s(M) = \xi_s^*(M)(-S)$ for any line bundle $M$ on $X_s$, of degree $\frak n$.   
Recall that $\xi_s$ is a composition $\xi_s = \pi_s\circ \varphi_s$
\[
\xymatrix { Y_r \ar@{->}[r]_{\varphi_s} \ar@/^0.5cm/[rr]^{\xi_s}  &  Y_r^* \ar@{->}[r]_{\pi_s} & X_s}
\]
where $\pi_s$ is the natural projection of the fiber product $Y_r^* = Y\times_X X_s$ and $\varphi_s$ is the normalization of $Y_r^*$, see Proposition~\ref{prop:etale} and Figure~\ref{fig:especcurves}.

The map $\Phi$ is defined using two correspondences, the former is a pullback and the second is an elementary transformation  followed  by a suitable reduction to the curve $Y_r$ (c.f. Section~\ref{sec:gens}).  Thus we need to follow the behavior of $M$ under these correspondences:
\[
\xymatrix { {\rm Pic}^{\frak n}(X_s) \ar@{->}[r]_{\pi_s^*} \ar@/^0.5cm/[rrr]^{\Phi_s}  &  {\rm Pic}^{2\frak n}(Y_r^*) \ar@{->}[r]_{elem_R} & {\rm Pic}^{2\frak n}(Y_r^*) \ar@{->}[r]_{\varphi_s^*} & {\rm Pic}^{\tilde{\frak n}}(Y_r)} .
\]

Let $(E, { l},\theta)$ be the Higgs bundle which corresponds to $M\in {\rm Pic}^{\frak n}(X_s)$. By Proposition \ref{prop:pullback}, $\pi_s^*M$ is the line bundle in ${\rm Pic}^{2{\frak n}}(Y_r^*)$  corresponding to  $(\pi^*(E, { l}),\pi^*\theta)$, and   it follows from Proposition \ref{twistM} that  the elementary transformation $\elem_R$ transforms  $\pi_s^*M$  into $\pi_s^*M\otimes \mathcal I$, where $\mathcal I$ is the ideal sheaf of $ (\zeta_r^*R)_{red}$ given by the 2:1 cover $\zeta_r: Y_r^*\to Y$.

The last operation is given by the normalization of $Y_r^*$. 
Since $\xi_s = \pi_s\circ\varphi_s$, then 
\[
\Phi_s(M) = \varphi_s^*(\pi_s^*M\otimes \mathcal I) = \xi_s^*(M)\otimes \varphi_s^*(\mathcal I) .
\]
The conclusion follows from the fact that $Y_r$ is smooth,  and then $\varphi_s^*(\mathcal I) = \mathcal O_{Y_r}(-S)$ is the ideal sheaf of $S=q_r^*R$, where $q_r: Y_r \to Y$ is the 2:1 cover.

To see that $\Phi$ has degree two onto its image,  first remark that  $\pi^*$ of (\ref{diag}) is a one-to-one map between Hitchin basis and   Proposition~\ref{prop:etalegeral}. 

\endproof

We point out that it would be interesting to compare the techniques used in this paper with the methods of \cite{Z22}.

We close the paper by presenting two situations where the modular map $\Phi$ is dominant, which deal with pullback of $SL_2$-Higgs bundles over $\P^1$ via $\pi: Y\to \P^1$, where $g_Y \in \{1, 2\}$.  By $SL_2$ we mean that $\theta: E\to E\otimes \omega_X(D)$ has vanishing trace and $E$ has trivial determinant line bundle.

 \subsection{The case $(g_X, g_Y)=(0, 1)$}\label{spechyper}
Let us consider the degree two elliptic cover $\pi: Y \to \P^1$, branched over four distinct points $\{0,1,\lambda, \infty\}$. Let us fix $t_1, t_2\in Y$ with  $\pi(t_1)=\pi(t_2)=t$, $t\notin\{0,1,\lambda, \infty\}$. 

On the  $\P^1$ side, we consider the moduli space $\mathcal H(\P^1, \Lambda)$ of irreducible $SL_2$ parabolic Higgs bundles over $(\P^1, \Lambda)$, where $\Lambda=0 + 1 + \lambda + \infty + t$. On the elliptic side, we consider the moduli  space $\mathcal H(Y, D)$ of irreducible $SL_2$ parabolic Higgs bundles over $(Y, D)$, where $D = t_1+t_2$.

It turns out that $\mathcal H(\P^1, \Lambda)$ and $\mathcal H(Y, D)$ have the same dimension: four. We can associate to a Higgs bundle in $\mathcal H(\P^1, \Lambda)$ a Higgs bundle in $ \mathcal H(Y, D)$ by doing the pullback $\pi^*$ followed by a composition of an elementary transformation $elem_R$ over the divisor \[
R=w_0+w_1+w_{\lambda}+w_{\infty}
\] 
formed by ramification points  of $\pi$.  This is our map $\Phi$. Note that, starting with a vector bundle $E$ of degree $0$ on $\P^1$, the transformed vector bundle on $Y$ has determinant $\mathcal O_Y(-R)$, then degree $-4$. Since $R\sim 4w_{\infty}$, we then perform a last transformation,  twisting  by $\mathcal O_Y(2w_{\infty})$, which results in trivial determinant.

This correspondence gives a map 
\begin{eqnarray*}
\Phi_0:\mathcal H(\P^1, \Lambda) \dashrightarrow \mathcal H(Y, D)
\end{eqnarray*}
which consists in $[\otimes\mathcal O_Y(2w_{\infty})]\circ\Phi$, and fits in  the following  commutative diagram 
 \[
\xymatrix { 
 \mathcal H(\P^1, \Lambda) \ar@{->}[d]_{\det} \ar@{-->}[r]^{\Phi_0}  &  \mathcal H(Y, D) \ar@{->}[d]^{\det} \\
\Gamma(\omega_{\P^1}^{\otimes 2}(\Lambda)) \ar@{->}[r]^{\pi^*}    &        \Gamma(\omega_{Y}^{\otimes 2}(t_1+t_2))
}
\]
(both spaces of quadratic differentials are two dimensional).

Over $\P^1$, the smooth spectral curve $X_s$ has genus $2$ and is branched over $6$ distinct points  $0,1,\lambda, \infty, t, \rho$. 
 The corresponding Hitchin fiber $\det^{-1}(s)$ is isomorphic ${\rm Pic}^{3}(X_s)$.  Over $Y$, the generic spectral curve $Y_r$ is a hyperelliptic  curve of genus $3$ branched over $4$ distinct points  $t_1, t_2, u_1, u_2$, with $\rho=\pi(u_1)=\pi(u_2)$. 
The fiber ${\rm det}^{-1}(r)$ over  $Y_r$  is an Abelian variety isomorphic to the Prym variety
\[
{\rm Prym}(Y_r/Y) = \left\{ M\in {\rm Pic}^{2}(Y_r)\;:\;\; \det((q_r)_*M)=\mathcal O_Y\right\} .
\]
Since $q_r$ is a ramified covering, then  ${\rm Prym}(Y_r/Y)$ is an irreducible variety, this follows from \cite[(iv) in p. 329]{Mu}.

We can choose local coordinates where 
\begin{displaymath}
\left\{ \begin{array}{ll}
X_s = \{(x,w)\quad : \quad w^2 = x(x-1)(x-\lambda)(x-t)(x-\rho)  \}\\
Y_r = \{ (x,y,z)\quad : \quad y^2 = x(x-1)(x-\lambda) \;\;\text{and}\;\; z^2 = (x-t)(x-\rho)\}
\end{array} \right. 
\end{displaymath}
and both curves fit in the commutative diagram below
\begin{eqnarray}
\xymatrix{
 &   & Y_r \ar[lld]_{2:1} \ar[d]_{2:1} \ar[rrd]^{{\text{2:1 \'etale: $\xi_s$}}}  & &\\
Y \ar[rrd]_{2:1}&  &   \P^1_{\rho} \ar[d]_{2:1}  &  & X_s \ar[lld]^{2:1}\\
  &  & \mathbb P^1 & &\\
}
\end{eqnarray}
where 
\[
\P^1_{\rho} = \{ (x,z)\quad : \quad z^2=(x-t)(x-\rho)\}. 
\]
The \'etale map $\xi_s: Y_r \to X_s$ is the same of Proposition~\ref{prop:etale},  explicitly given here by $(x,y,z)\mapsto(x,yz)$.

In this context, Theorem~\ref{thm:main} can be written as

\begin{cor}\label{cor:P1C}
The rational map  $\Phi_0:\mathcal H(\P^1, \Lambda) \dashrightarrow \mathcal H(Y, D)$ has degree two. Moreover, its restriction to a general fiber of the Hitchin map is a morphism of degree two given by 
\begin{eqnarray*}
\Phi_{0,s}: {\rm Pic}^3(X_s) &\to& {\rm Prym}(Y_r/Y)\\ 
             M &\mapsto& \xi_s^*(M) (S_0)
\end{eqnarray*} 
where $S_0 = q_r^*(2w_{\infty} - R)$.  
\end{cor}

\proof
As mentioned above, the map $\Phi_0$ consists of $[\otimes\mathcal O_Y(2w_{\infty})]\circ\Phi$ where $\Phi$ is the map of Theorem~\ref{thm:main}. Let us fix a point $(E,l,\theta)$ in the image of $\Phi_0$ with characteristic coefficient $r$. By construction, we have $\det E = \mathcal O_Y$, then it follows by BNR correspondence that $(E,l,\theta)$ corresponds to a line bundle $N\in {\rm Pic}^2(Y_r)$ satisfying $\det ((q_r)_*N) = \det E = \mathcal O_Y$. Hence, $N\in {\rm Prym} (Y_r/Y)$ and the proof of the corollary follows from Theorem~\ref{thm:main}.
\endproof

In a forthcoming paper, we will give the full description
(at singular Hitchin fibers) of the map $\Phi_0$.

\subsection{The case $(g_X, g_Y) = (0, 2)$}

We now take a degree two morphism $\pi: Y\to \P^1$ from a curve $Y$ of genus $2$. It is branched over $6$ distinct points $t_1, \dots, t_6$ of $\P^1$.   Let $B= t_1 + \cdots + t_6$ be the reduced divisor on $\P^1$ defined by those points, and let $R = w_1+\cdots + w_6$ be the divisor on $Y$ formed by ramification points of $\pi$. 

Let us consider the moduli space $\mathcal H(\P^1, B)$ of irreducible $SL_2$-Higgs bundles over $(\P^1, B)$. Over $Y$, we consider the moduli space $\mathcal H(Y)$ of irreducible holomorphic $SL_2$-Higgs bundles over $Y$, without parabolic points. These moduli spaces have both dimension $6$. 

The spectral curve on the $\P^1$ side has genus $3$ and the corresponding smooth  Hitchin fiber is isomorphic to ${\rm Pic}^4(X_s)$. Over $Y$, the spectral curve $Y_r$ has genus $5$ and the smooth Hitchin fiber $\mathcal H(Y)$                                                                                                                                                                is isomorphic to  
\[
{\rm Prym}(Y_r/Y) = \left\{ M\in {\rm Pic}^{2}(Y_r)\;:\;\; \det((q_r)_*M)=\mathcal O_Y\right\} 
\]
which is three dimensional, and irreducible (\cite[(iv) in p. 329]{Mu}).

We consider the correspondence $\Phi_0: \mathcal H(\P^1, B) \dashrightarrow \mathcal H(Y)$ given by $[\otimes L_0]\circ elem_R\circ \pi^*$, where $L_0$ is a square root of $\mathcal O_Y(R)$. Note that the following diagram is commutative   
\[
\xymatrix { 
 \mathcal H(\P^1, B) \ar@{->}[d]_{\det} \ar@{-->}[r]^{\Phi_0}  &  \mathcal H(Y) \ar@{->}[d]^{\det} \\
\Gamma(\omega_{\P^1}^{\otimes 2}(B)) \ar@{->}[r]^{\pi^*}    &        \Gamma(\omega_{Y}^{\otimes 2})
}
\]
where the Hitchin basis are three dimensional. From Theorem~\ref{thm:main}, we obtain a result similar  to Corollary~\ref{cor:P1C}:
\begin{cor}\label{cor:g=2}
The rational map  $\Phi_0:\mathcal H(\P^1, B) \dashrightarrow \mathcal H(Y)$ has degree two and its restriction to a general fiber of the Hitchin map consists of  
\begin{eqnarray*}
\Phi_{0,s}: {\rm Pic}^4(X_s) &\to& {\rm Prym}(Y_r/Y)\\ 
             M &\mapsto& \xi_s^*(M) \otimes q_r^*(L_0(-R)) .
\end{eqnarray*} 
\end{cor}

\proof
Similar to the proof of Corollary~\ref{cor:P1C}. 
\endproof

\end{document}